\documentclass[12pt]{article}
\usepackage{authblk}
\usepackage[utf8]{inputenc}
\usepackage[cal=boondoxo]{mathalfa}
\usepackage{amssymb}
\usepackage{comment}
\usepackage{amsmath}
\usepackage{graphicx}
\usepackage{mathtools}
\usepackage{wrapfig}
\usepackage{xcolor}
\usepackage{verbatim}
\usepackage{indentfirst}
\usepackage{enumerate}
\usepackage{a4}
\usepackage{graphics}
\usepackage{epsfig}
\usepackage[colorlinks=true,linkcolor=black,citecolor=black,urlcolor=black]{hyperref}
\newtheorem{thm}{Theorem}[section]
\newtheorem{lm}[thm]{Lemma}

\numberwithin{equation}{section}

\newcommand{\R}{\mathbb{R}}
\newcommand{\D}{\mathbb{D}}
\newcommand{\T}{\mathbb{T}}

\newcommand{\N}{\mathbb{N}}
\newcommand{\El}{\mathcal{E}}

\newcommand{\norm}[1]{\Vert #1\Vert}

\newcommand{\Dn}[1]{\mathbb{D}^{#1}}
\newcommand{\Cn}[1]{\mathbb{C}^{#1}}
\newcommand{\Bn}[1]{\mathbb{B}^{#1}}

\newcommand{\br}[1]{\{#1\}}

\newcommand{\hols}[2]{\mathcal{O}(#1,#2)}
\newcommand{\address}[1]{\par\noindent\textbf{Address: }#1}
\newcommand{\email}[1]{\par\noindent\textbf{Email: }\texttt{#1}}

\def\iint{\operatorname{int}}

\begin{document}
\thispagestyle{plain}
\begin{center}
    \Large
    \textbf{The uniqueness problem for the 3-Point Nevanlinna-Pick extremal interpolation problem in the unit Euclidean ball $\Bn{d}$}

    \vspace{0.4cm}
    \text{Dariusz Piekarz}\footnote{Partially supported by the National Center of Science, Poland, Preludium Bis 2 grant no. 2020/39/O/ST1/00866}
\end{center}       
 \vspace{0.9cm}
 
\begin{abstract}Agler and McCarthy studied the uniqueness of a 3-point interpolation problem in the bidisc. This note aims to solve an analogous problem in the unit Euclidean ball in an arbitrary dimension.
\end{abstract}

\vspace{0.4cm}
\section{Introduction}
In this paper, we study the uniqueness problem for an extremal 3-point Nevanlinna-Pick interpolation problem in the unit Euclidean ball $\Bn{d}$. A similar problem was addressed by Agler and McCarthy for the 3-point Nevanlinna-Pick interpolation problem in the bidisc.

Let us recall that an $N$-point Nevanlinna-Pick interpolation problem is a classical question that can be formulated as follows:
Given a bounded domain $\Omega$ in $\Cn{d},$ consider $N$ pairwise distinct points $z_{i} \in \Omega$ and numbers $\zeta_{i} \in \D$ (not necessarily pairwise distinct). Determine whether there exists a holomorphic function $f \in \hols{\Omega}{\D}$ that interpolates $\Omega \ni z_{i} \mapsto \zeta_i \in \D$ for all $i = 1, \dots, N.$

An interpolation problem is called \textit{extremal} if a solution exists, but there is no solution whose image lies relatively compactly on the unit disc $\D$. 

Let us introduce some notation. An open disc at $z_0$ with a radius $r$ is denoted by $\D(z_0,r).$ A particular case of the unit disc in the complex plain will be denoted by the symbol $\D,$ while $\D_{*}$ denotes the punctured disc, i.e., $\D_{*} := \D \setminus \br{0}.$ The unit circle is denoted by $\T.$ Furthermore, $\Bn{d}$ is the complex unit Euclidean ball in $\Cn{d}.$ A topological interior of a set $A$ is denoted as $\iint(A),$ and the space of holomorphic functions mapping from $\Omega_1 $ to $\Omega_2$ is denoted by $\hols{\Omega_1}{\Omega_2}.$

\section{The results}

\par An interpolation problem $\mathbb D\to \mathbb D$ was studied by Pick in 1916 and independently by Nevanlinna in 1919 who also focused on the uniqueness part. Today, this problem is well understood and can be proven by a simple induction, e.g. using Schur's algorithm. However, this approach cannot be generalized to higher-dimensional settings. The only far-reaching result is known only for the bidisc. Actually, in \cite[Chapter 8.3]{-1}, Agler extended Pick’s theorem to the space $H^{\infty}(\Dn{2})$, the bounded analytic functions on the bidisc. Based on this result, Agler and McCarthy were able to understand the uniqueness of an extremal solution to the $3$-point Nevanlinna-Pick interpolation problem for the bidisc. This statement was later refined by Kosiński \cite{2} who used a particular solution of the $3$-point Nevanlinna-Pick interpolation problem for the polydisc $\Dn{d}$ in terms of generalized geodesics. In 2018, Kosiński and Zwonek, in \cite{3}, presented an analogous solution to the $3$-point Nevanlinna-Pick interpolation problem in the unit Euclidean ball $\Bn{d}$. Precisely, they proved the following
\begin{thm}[Kosiński, Zwonek]\label{kzw:1}
If the $3$-point Pick interpolation problem
\[
\Bn{d} \to \D, \quad w_j \mapsto \lambda_j, \quad j = 1, 2, 3,
\]
is extremal, then, up to a composition with automorphisms of $\Bn{d}$ and $\D$, it is interpolated by a function belonging to one of the classes:
\[
\mathcal{F}_{D} = \left\{ (z_1, \dots, z_d) \mapsto \frac{2z_1(1 - \tau z_1) - \overline{\tau} \omega^2 z_2^2}{2(1 - \tau z_1) - \omega^2 z_2^2} : |\tau| = 1, |\omega| \leq 1 \right\},
\]
\[
\mathcal{F}_{ND} = \left\{ (z_1, \dots, z_d) \mapsto \frac{z_1^2 + 2\sqrt{1 - a^2} z_2}{2 - a^2} : a \in [0, 1) \right\}.
\]
Therefore, the solutions have the form $m \circ f \circ \phi$, where $m \in \text{Aut}(\D)$, $\phi \in \text{Aut}(\Bn{d})$, and $f \in \mathcal{F}_D$ or $f \in \mathcal{F}_{ND}$.
\end{thm}
Indices of the classes in the theorem above are abbreviations for "Degenerate" and "Non-Degenerate". In the case of the $3-$point extremal problem, one can distinguish two types of the problem: if any 2-point subproblem mapping the pair $(z_i, z_j)$ to $(\lambda_i, \lambda_j)$ is extremal, we call the problem \textit{degenerate}. Otherwise, it is called \textit{non-degenerate.}

In fact, degeneracy implies that for a bounded and convex domain $\Omega$, the Carath\'eodory distance  
\[
c_\Omega(z, w) := \sup \left\{ \rho(F(z), F(w)) : F \in \hols{\Omega}{\mathbb{D}}\right\},
\]
where $\rho$ is the Poincar\'e distance in $\mathbb{D}$, satisfies the following condition for a $2$-point extremal subproblem $(z_1, z_2)$ mapped to $(\lambda_1, \lambda_2)$:  
\[
c_\Omega(z_1, z_2) = c_\mathbb{D}(\lambda_1, \lambda_2).
\]

In the case of $\Omega = \Bn{d}$, as shown later, one can additionally apply automorphisms of the unit Euclidean ball to assume that the extremal subproblem has the form
\[
\begin{cases}
\Bn{d}\ni 0\mapsto 0\in\D \\
\Bn{d}\ni z \mapsto \sigma\in\D, 
\end{cases}
\]
In this particular case, it means that the degeneracy can be described by the condition $\norm{z}_2=|\sigma|.$ 
\par It is natural to ask whether the solution to the extremal problem is unique. For the case of bidisc $\Dn{2}$ the answer was found by Agler and McCarthy in \cite{0}. They proved that in the non-degenerate case the extremal problem always has a unique solution. Kosiński extended the Agler and McCarthy results to polydisc $\Dn{d},d\geq3$ proving that in the non-degenerate case the extremal problem never has a unique solution.

\par On a separate note, it should be highlighted that the Nevanlinna-Pick interpolation problem has intriguing applications in other areas of mathematics, such as the von Neumann inequality (see, for instance, applications of Kosiński's and Kosiński-Zwonek's results in \cite{8}, \cite{4}) or extensions properties (ex. Kosiński and McCarthy, \cite{7}).
\par In this paper, we investigate the uniqueness of solutions to the 3-point Nevanlinna-Pick extremal interpolation problem in the unit Euclidean ball $\Bn{d}$. The main result has a simple formulation and differs from the one for the bidisc:
\begin{thm} 
    Solutions of the 3-point Nevanlinna-Pick extremal interpolation problem in the Euclidean ball are never unique. 
\end{thm}

\section{Theory of the holomorphic functions of the unit Euclidean ball}
\subsection{Automorphisms of the unit Euclidean ball $\Bn{d}$}
 Let us begin recalling the form of the automorphism group of the unit disc $\D$ and the unit Euclidean ball $\Bn{d}$.
The group of automorphisms of the unit disc is given by 
$$\text{Aut} (\mathbb{D}) = \left\{ \mathbb{D} \ni z \mapsto \omega \frac{z - a}{1 - \overline{a}z} \in \mathbb{D} \, : \, \omega \in \T, \, a \in \mathbb{D} \right\};$$ while the group of automorphisms of the unit Euclidean ball $\Bn{d}$ is described as
\[
 \ \text{Aut}(\mathbb{B}^d) = \left\{ U \circ \phi_a : U \in \mathcal{U}(\mathbb{C}^d), \ a \in \mathbb{B}^d \right\},\] \[ \text{Aut}_0 (\mathbb{B}^d) = \mathcal{U}(\mathbb{C}^d),
\]
where for $a \in \Bn{d}, \ \phi_a(z)$ is defined as:

\begin{align}
\phi_a(z) := 
\begin{cases} 
\frac{1}{\|a\|^2} \frac{\sqrt{1 - \|a\|^2} \big( \|a\|^2 z - \langle z, a \rangle a \big) - \|a\|^2 a + \langle z, a \rangle a}{1 - \langle z, a \rangle}, & \text{if } a \neq 0, \\[8pt]
\mathrm{id}, & \text{if } a = 0.
\end{cases}
\end{align}
Here $\langle \cdot, - \rangle$ denotes the standard scalar product in $\mathbb{C}^d$ and  $\mathcal{U}(\mathbb{C}^d$)  is the group of unitary automorphisms of $\mathbb{C}^d$.

The automorphism group $\text{Aut}(\Bn{d})$ is generated by finite composition of the unitary mapping with the automorphisms $\phi_{(a_1,0,...,0)},$ i.e., with mappings of the form:
$$\phi_{(a_1,0,...,0)}(z_1,...,z_d)=\left(\frac{z_1-a_1}{1-\overline{a_1}z_1},\sqrt{1-|a_1|^2}\frac{z_2}{1-\overline{a_1}z_1},...,\sqrt{1-|a_1|^2}\frac{z_{d}}{1-\overline{a_1}z_1}\right).$$
\subsection{Extremality of the 3-point Nevanlinna-Pick interpolation problem}
Suppose that the N-point Nevanlinna-Pick interpolation problem has the following form 
\[
\begin{cases}
    \Omega\ni z_1 \mapsto\zeta_1\in\D\\
    ...\\
    \Omega\ni z_N \mapsto\zeta_N\in\D,
\end{cases}
\]
with $N \geq 2.$ Additionally, assume that at least two points from $\zeta_1,...,\zeta_N$ are distinct (otherwise, the solution to the problem is a constant function).
Suppose that this problem can be solved by some holomorphic function $f$. The problem is not necessarily extremal, we can from this construct an extremal problem and find its solution. For this consider modified problem 
\begin{align}\label{modpr}
\begin{cases}
    \Omega\ni z_1 \mapsto t\zeta_1\in\D\\
    ...\\
    \Omega\ni z_N \mapsto t\zeta_N\in\D,
\end{cases}
\end{align}
where $t>0.$
Consider a set $\tau:=\br{t>0: \eqref{modpr} \ \text{has a solution}}.$ Then $\tau\neq\emptyset$ (as $1\in\tau$) and $\tau$ is bounded as a range of any solution must be contained in the unit disc $\D$. Consider a sequence of $\br{t_n}_{n\in\N}\subset\tau$ and $t_n\rightarrow t_0:=\sup\tau$ as $n\rightarrow+\infty.$ Then to each $t_n$ one has a corresponding solution $F_n\in\hols{\Omega}{\D}$ of \eqref{modpr}. Therefore, we have uniformly bounded sequence of holomorphic functions $\br{F_n}_{n\in\N}.$ Due to the Montel theorem there exists a subsequence $\br{t_{n_k}}_{k\in\N}\subset\br{t_{n}}_{n\in\N}$ such that it converges uniformly to a holomorphic function $F:\Omega\rightarrow\overline{\D}$ i.e., $F_{n_k}\rightarrow F.$ As numbers $\zeta_1,...,\zeta_N$ are taken from the unit disc $\D,$ and at least two from these points are distinct the Maximum Modulus Principle asserts that in fact $F$ is not constant and $F(\Omega)\subset\D.$ The function $F$ satisfies $F(z_1)=t_0\zeta_1,...,F(z_N)=t_0\zeta_N$ and the problem 
\begin{align*}\label{modpr}
\begin{cases}
    \Omega\ni z_1 \mapsto t_0\zeta_1\in\D\\
    ...\\
    \Omega\ni z_N \mapsto t_0\zeta_N\in\D
\end{cases}
\end{align*}
is extremal.
The above procedure can be visualized as in the pictures below ($N=3$):
\begin{figure}[h!]
    \centering
    \includegraphics[width=0.8\textwidth]{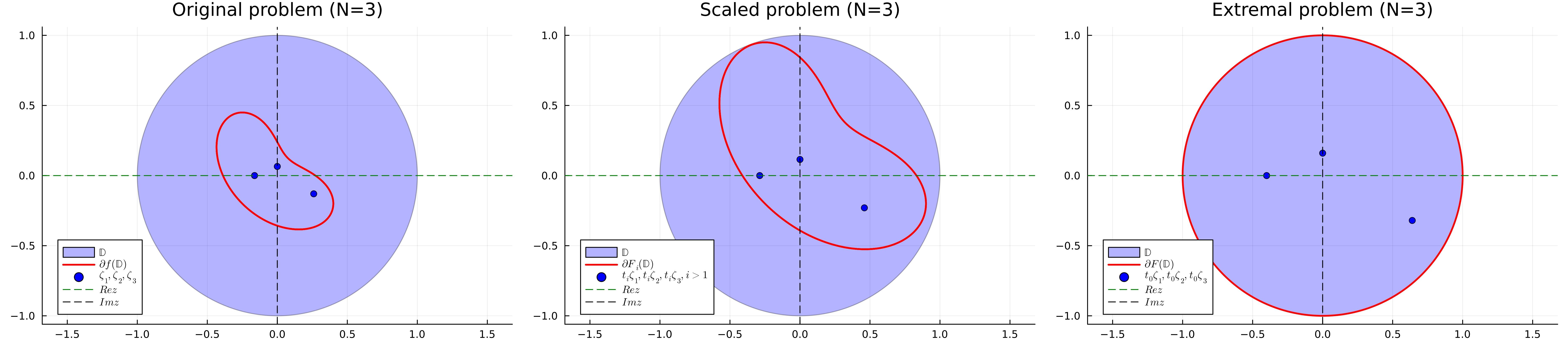}
\end{figure}

As can be seen in the above picture, the solutions to the extremal problem have their images dense in $\D.$ This can be shown in the spirit of \cite[Lemma 9.4]{-2}.
\subsection{3-point Nevanlinna-Pick problem in \texorpdfstring{$\Bn{d}$}{} and its reductions}
\par First, let us state the problem for the unit Euclidean ball $\Bn{d}$. Consider the following interpolation conditions:
\[
\begin{cases} \Bn{d}\ni v \mapsto \chi\in\D \\ \Bn{d}\ni w \mapsto \kappa\in\D \\ \Bn{d}\ni z \mapsto \sigma\in\D, \end{cases}\] where $z, w, v$ are pairwise distinct and $\chi, \kappa, \sigma\in\D$ (not necessarily pairwise distinct). Applying the automorphisms of the unit ball and the unit disc, we can assume without loss of generality that $v = 0$ and $\chi = 0$. Additionally, we can assume that at least one of $ \sigma, \kappa$ is different from$0,$ otherwise the solution would be a constant function, which is not an interesting case. Using rotations, which are automorphisms of the ball, we first rotate the point $z$ so that it lies in the plane spanned by the coordinate vectors $e_1 := (1, 0, \dots, 0)$ and $e_2 := (0, 1, 0, \dots, 0)$. After this rotation, the point $z$ takes the form $(z_1, z_2, 0, \dots, 0)$. Note that the point $w$ is also rotated to another point within $\Bn{d}$. The problem has the following form: 
\[
\begin{cases}
\Bn{d}\ni (0,0,0,...,0)  \mapsto 0\in\D \\
\Bn{d}\ni (w_1,w_2,...,w_d) \mapsto \kappa\in\D  \\
\Bn{d}\ni (z_1,z_2,0,...,0) \mapsto \sigma\in\D, 
\end{cases}
\]
Then, within the plane, we apply a rotation so that the point $(z_1, z_2, 0, \dots, 0)$ is transferred to a point on the axis determined by $e_1$. Note that the point $0$ remains invariant under these rotations, and the point $w$ is mapped to another point.

Thus, we can further assume that $z = (z_1, 0, \dots, 0)$ with $z_1 > 0$.The problem is now modified to:
\[
\begin{cases}
\Bn{d}\ni (0,0,0,...,0) \mapsto 0\in\D \\
\Bn{d}\ni (w_1,w_2,...,w_d) \mapsto \kappa\in\D  \\
\Bn{d}\ni (z_1,0,0,...,0) \mapsto \sigma\in\D, 
\end{cases}
\]
Applying rotations around the axis determined by the vector $e_1$, we can also map the point $w$ into the plane spanned by $e_1$ and $e_2$. Hence, we can assume that $w$ has the form $w = (w_1, w_2, 0, \dots, 0)$.

Therefore, the problem has the following form:
\[
\begin{cases}
\Bn{d}\ni (0,0,...,0)\mapsto 0\in\D \\
\Bn{d}\ni (w_1,w_2,0,...,0) \mapsto \kappa\in\D  \\
\Bn{d}\ni (z_1,0,0,...,0) \mapsto \sigma\in\D, 
\end{cases}
\]

Henceforth, one can assume that $d = 2$. If $F: \Bn{2} \rightarrow \D$ is an extremal solution in the 2-dimensional ball, then $$\Bn{d} \ni (z_1, \dots, z_d) \mapsto F(z_1, z_2) \in \D$$ is a solution in $\Bn{d}$ for $d > 2$. Conversely, if $F(z_1,...,z_d)$ is a solution to the above problem in $\Bn{d},d>2,$  then $G(z_1,z_2):=F(z_1,z_2,0,...,0)$ is the solution in $\Bn{2}.$

Next, consider a solution $F$ to the extremal problem satisfying the following conditions:
\[
\begin{cases}
\Bn{2}\ni (0,0) \mapsto 0\in\D \\
\Bn{2}\ni (w_1,w_2) \mapsto \kappa\in\D  \\
\Bn{2}\ni (z_1, 0) \mapsto \sigma\in\D, 
\end{cases}
\]
In the non-degenerate case described above, these reductions are sufficient. However, in the degenerate case, the problem can be simplified even further.
Due to the degeneracy of the problem, it follows that $z_1 = \norm{(z_1, 0)}_2 = |\sigma|$, as $z_1 > 0$. Applying a rotation of the unit disc, we can assume that $\sigma = z_1$ (this also affects the value of $\kappa$ as well).

Now, consider the function $f: \D \ni \lambda \mapsto F(\lambda, 0) \in \D$. For this function, $f(0) = F(0, 0) = 0$ and $f(z_1) = F(z_1, 0) = z_1$. Since $f$ is a holomorphic self-mapping of the unit disc $\D$, applying the Schwarz lemma implies that $f(\lambda) = F(\lambda, 0) = \lambda$.

This result shows that the third condition in the considered interpolation problem, $(z_1, 0) \mapsto z_1$, can be replaced with any condition $(\lambda_0, 0) \mapsto \lambda_0$, where $\lambda_0 \in \D$.

Thus, the original problem can be equivalently reformulated as follows:
\[
\begin{cases}
\Bn{2}\ni (0, 0) \mapsto 0\in\D \\
\Bn{2}\ni (w_1,w_2) \mapsto \kappa\in\D  \\
\Bn{2}\ni (w_1, 0) \mapsto w_1\in\D, 
\end{cases}
\]
One further reduction is possible. Using an automorphism of the unit ball, $\text{Aut}(\Bn{2})$, of the form  $$\phi_{(w_1,0)}(z_1,z_2)=\left(\frac{z_1-w_1}{1-\overline{w_1}z_1},\sqrt{1-|w_1|^2}\frac{z_2}{1-\overline{w_1}z_1}\right)$$ and applying an appropriate automorphism of the unit disc $\D$, the problem can be further reduced to:
\begin{align}\label{problem}
\begin{cases}
    \Bn{2}\ni (0, 0) \mapsto 0\in\D \\
    \Bn{2}\ni (0, w_2) \mapsto \kappa\in\D \\
    \Bn{2}\ni (z_1, 0) \mapsto z_1\in\D. 
\end{cases}    
\end{align}

\section{The proof}

To prove the main result, we shall consider two cases depending on whether the problem is degenerate or not. 

\subsection{Non-degenerate case}
In this section, our aim is to prove that if the problem
\[
\begin{cases}
\Bn{2}\ni (0,0) \mapsto 0\in\D \\
\Bn{2}\ni (w_1,w_2) \mapsto \kappa\in\D  \\
\Bn{2}\ni (z_1, 0) \mapsto \sigma\in\D, 
\end{cases}
\]
is extremal and non-degenerate, then its solution is not unique.
\par Before we start, recall that an analytic disc $f : \mathbb{D} \to \Omega$ is called a \emph{complex $N$-geodesic} if there exists a holomorphic function $F:\Omega\rightarrow\D$ such that $b := F \circ f$ is a non-constant Blaschke product of degree at most $N - 1$.

The nonuniqueness in this case can be deduced from Kosiński and Zwonek in \cite{5} and \cite{3}. We will follow their approach. In the non-degenerate case, the points $(0, 0)$, $(z_1, 0)$, and $(w_1, w_2)$ lie in the range of the complex 3-geodesic \begin{align}\label{3geo}\varphi:\D \ni \lambda \mapsto \left(a\lambda, \sqrt{1-a^2}\lambda^2\right) \in \Bn{2}\end{align} for some $a \in [0,1)$. Its left inverse is given by \begin{align}\label{ngleft} F_a(z, w) = \frac{1}{2-a^2}\left(z^2 + 2\sqrt{1-a^2}w\right).\end{align}

We will prove that for any $\varepsilon\in[0,1)$ the function $$F_{a,\varepsilon}(z,w):=\frac{F_a(z,w)}{\sqrt{1-\frac{\varepsilon^2}{(2-a^2)^2}(a^2w-bz^2)^2}}$$ is a left inverse to the 3-geodesics \eqref{3geo}, where $b:=\sqrt{1-a^2}$.
First we prove that the range of $F_{a,\varepsilon}$ lies in the unit disc $\D$ for $\varepsilon\in(0,1).$
For this we must prove the following inequality
$$|z^2+2bw|^2+\varepsilon^2|a^2w-bz^2|^2<(2-a^2)^2.$$
The proof of the inequality above can be reduced to the proof of the following inequality:
$$f(r)\leq(2-a^2)^2=f(b), r\in[0,1],$$ where $f(r):=(1-r^2+2br)^2+\varepsilon^2(a^2r-b(1-r^2))^2.$ Indeed, due to Maximum Modulus Principle one can assume that $|z|^2+|w|^2=1.$ Set $r:=|w|,$ and $z^2=(1-r^2)e^{i\theta},\theta\in(0,\pi],$ $w=re^{i\alpha},\alpha\in(0,2\pi]$. 
Then we calculate that 
\begin{align*}
|z^2+2bw|^2+\varepsilon^2|a^2w-bz^2|^2&=(z^2+2bw)(\overline{z}^2+2b\overline{w})+\varepsilon^2(a^2w-bz^2)(a^2\overline{w}-b\overline{z}^2)\\&=(1-r^2)+4br(1-r^2) Re (e^{i(\theta-\alpha)})+4b^2r^2 \\ &+\varepsilon^2(a^4r^2-2a^2br(1-r^2)\cos{(\theta-\alpha)}+b^2(1-r^2)^2)\\ &=(1+\varepsilon^2b^2)(1-r^2)^2+(4b^2+\varepsilon^2a^4)r^2\\ &+2b(2-a^2\varepsilon^2)r(1-r^2)\cos{(\theta-\alpha)}.
\end{align*}
Provided that $2-a^2\varepsilon^2>0$ one can estimate the last expression from above by 
$$(1+\varepsilon^2b^2)(1-r^2)^2+(4b^2+\varepsilon^2a^4)r^2+2b(2-a^2\varepsilon^2)r(1-r^2)=f(r).$$

Hence, $|z^2+2bw|^2+\varepsilon^2|a^2w-bz^2|^2\leq f(r)$, and if the inequality for $f$ is true, the original inequality is true as well. The polynomial $f$ has degree 4, and has a positive leading coefficient. Furthermore, $f(r)\geq0$ for any $r\in\R$. The derivative of $f$ is
$$f'(r)=2(b-r)(-2(b^2\varepsilon^2+1)r^2+b(4-3\varepsilon^2+b^2\varepsilon^2)r + 2-\varepsilon^2+b^2\varepsilon^2).$$ We prove that the quadratic function $$q(r):=-2(b^2\varepsilon^2+1)r^2+b(4-3\varepsilon^2+b^2\varepsilon^2)r + 2-\varepsilon^2+b^2\varepsilon^2,$$ has no roots in the interval $[0,1).$  For this we use Viet\'e's formula for products of roots.
First, we observe that the discriminant of the square polynomial is greater than zero:
$$b^2(4-3\varepsilon^2+b^2\varepsilon^2)^2+8(b^2\varepsilon^2+1)(2-\varepsilon^2+b^2\varepsilon^2)>0,$$ then Viet\'e's formula for product of roots is 
$$r_1r_2=\frac{2-\varepsilon^2+b^2\varepsilon^2}{-2(b^2\varepsilon^2+1)}<0,$$ as the nominator is positive and the denominator is negative. The additional condition to make sure that the roots of $q$ are outside the interval $[0,1)$ is to check that $q(1)\geq0.$ We have \begin{align*}
q(1)&=4b+\varepsilon^2(b^3-b^2-3b-1)\geq4b\varepsilon^2+\varepsilon^2(b^3-b^2-3b-1)\\ &=\varepsilon^2(b^3-b^2+b-1)=\varepsilon^2(b-1)(b^2+1)\geq0.
\end{align*}

\par Therefore, it suffices to find the maximum for $f$ and consider the end of interval values of $f$ i.e., values at points $0,1$ as $f$ has potentially extreme points at $0,1$ and $b$. We have $f(0)=1+\varepsilon^2b^2$, $f(1)=4b^2+\varepsilon^2(1-b^2)^2$, $f(b)=(1+b^2)^2$. At point $b$ we have the extreme point, which is a maximum as $f''(b)=-4(2-a^2)+2\varepsilon^2(2-a^2)^2$. One has $f''(b)<0$ if and only if $\varepsilon^2<\frac{2}{2-a^2}$, which is satisfied if $\varepsilon\in(0,1).$

Therefore, to prove the inequality we must show that $$f(1)=4b^2+\varepsilon^2(1-b^2)^2\leq f(b)=(1+b^2)^2=(2-a^2)^2.$$ This inequality can be rearranged to $(1-\varepsilon^2)(b^2+1)^2+\varepsilon^2>0$, which is clearly true. Therefore, we can estimate the norm of $F_{a,\varepsilon}$ using the inequality we proved above $$|F_{a}(z,w)|^2+\frac{\varepsilon^2}{(2-a^2)^2}|a^2w-bz^2|^2<1.$$ Applying the triangle inequality
$$|F_{a}(z,w)|^2<\left|1-\frac{\varepsilon^2}{(2-a^2)^2}(a^2w-bz^2)^2\right|$$  and this is equivalent to $|F_{a,\varepsilon}(z,w)|<1.$ Therefore, $F_{a,\varepsilon}\in\hols{\Bn{2}}{\D}$ for any $\varepsilon\in(0,1).$ Clearly, one has $F_{a,\varepsilon}\circ\varphi=id_{\D},$ for any $\varepsilon\in(0,1).$ The nonuniqueness follows.

\subsection{Degenerate case}

In this section, we will show that the 3-point extremal and degenerate Nevanlinna-Pick interpolation problem never has a unique solution.
Recall that the problem can be reduced to the form
\begin{align}\label{deg_red}
\begin{cases}
\Bn{2}\ni (0, 0) \mapsto 0\in\D \\
\Bn{2}\ni (0, w_2) \mapsto \kappa\in\D  \\
\Bn{2}\ni (z_1, 0) \mapsto z_1\in\D.
\end{cases}
\end{align}

In \cite{3}, Kosiński and Zwonek characterized the degeneracy of the problem by describing the possible values of the set

\begin{align}\label{bw}
\mathcal{B}(w) := \{ F(w) : F \in \mathcal{O}(\mathbb{B}^2, \mathbb{\D}), F(z_1, 0) = z_1, z_1 \in \mathbb{\D} \}.
\end{align}
They showed in \cite{3} the following
\begin{thm}[Kosiński, Zwonek]\label{Bw_descr}
Let $w \in \mathbb{B}^2$. Then
\begin{equation}
\mathcal{B}(w) = m_{w_1} \left( \mathcal{B} \left( 0, \frac{w_2}{\sqrt{1 - |w_1|^2}} \right) \right)
= m_{w_1} \left( \overline{\D} \left( 0, \frac{|w_2|^2}{2 - 2|w_1|^2 - |w_2|^2} \right) \right).
\end{equation}
In particular, the set $\mathcal{B}(w)$ is a closed Euclidean disc.
Moreover, the extremal 3-point Pick interpolating functions in the degenerate case may be chosen from a nice class of functions. More precisely,
\[
\mathcal{B}(w) = \left\{ F_{\tau, \omega}(w)=\frac{2w_1(1 - \tau w_1) - \overline{\tau} \omega^2 w_2^2}{2(1 - \tau w_1) - \omega^2 w_2^2} : |\tau| = 1, |\omega| \leq 1  \right\},
\] where functions $F_{\tau, \omega}$ are described by the class $\mathcal{F}_D$ in Theorem \ref{kzw:1}.
 
\end{thm}
\par In view of the theorem above, the set $\mathcal{B}(w)$ is a closed Euclidean disc. We shall see that $$\partial\mathcal{B}(w)=\br{F_{\tau,\omega_2}(w): \tau\in\T},$$ where $\omega_2\in\T$ is such that $|w_2|^2=\omega_2w_2^2.$

For this, using a simplification from \eqref{deg_red}, it suffices to consider
$$\mathcal{B}(0,w_2) = \{ F(0,w_2) : F \in \mathcal{O}(\mathbb{B}^2, \mathbb{\D}), F(z_1, 0) = z_1, z_1 \in \mathbb{\D} \}.$$ 
Consider a function $\phi := \phi_{\tau,w_2}:\overline{\D} \ni \omega \mapsto F_{\tau,\omega}(0,w_2) \in \overline{\D}$, which is of the form $$\phi(\omega)=\frac{-\tau\omega^2w_2^2}{2-\omega^2w_2^2},\ \omega\in \overline{\mathbb D}.$$ The maximum of $\phi$ is attained on the circle $\T$. Moreover, for any $w_2$, there is a $\omega_2$ such that $|w_2| = \omega_2 w_2$, and then $$|\phi(\omega_2)|=\frac{|w_2|^2}{2-|w_2|^2}\in\partial\mathcal{B}(0,w_2).$$
Hence, if $\omega=\omega_2 \in \T$, then $F_{\tau,\omega}(w) \in \partial \mathcal{B}(0,w_2)$; otherwise, $F_{\tau,\omega}(w)$ lies in the open disc $\text{int}\mathcal{B}(0,w_2)$.

Let $\kappa\in\mathcal{B}(0,w_2),$ we consider two cases: $\omega=\omega_2$ and $\omega\neq\omega_2.$

\subsection{Case: $\omega\neq\omega_2$}

In view of the argument above, if $\omega\neq\omega_2$, then $\kappa$ lies in $\iint(\mathcal{B}(0,w_2))$.

In this situation, there exist $\mu, \nu \in \D$ such that $\mu \neq \nu$, $\mu \neq \kappa \neq \nu$, and $t \in (0, 1)$ such that $t\mu + (1-t)\nu = \kappa$, where $\mu, \nu \in \iint(\mathcal{B}(0,w_2))$.
Therefore, we can consider related problems:
\begin{align}\label{mpr}
\begin{cases}
\Bn{2}\ni (0, 0) \mapsto 0\in\D \\
\Bn{2}\ni (0, w_2) \mapsto \mu\in\D  \\
\Bn{2}\ni (z_1, 0) \mapsto z_1\in\D,
\end{cases}
\end{align}

and

\begin{align}\label{npr}
\begin{cases}
\Bn{2}\ni (0, 0) \mapsto 0\in\D \\
\Bn{2}\ni (0, w_2) \mapsto \nu\in\D  \\
\Bn{2}\ni (z_1, 0) \mapsto z_1\in\D.
\end{cases}
\end{align}
As $\mu, \nu \in \text{int}\mathcal{B}(0,w_2)$, both problems \eqref{mpr} and \eqref{npr} are extremal, degenerate, and extremally solvable. Denote their solutions by $F_\mu$ and $F_\nu$, correspondingly. Since $F_\mu(0, w_2) = \mu \neq \nu = F_\nu(0, w_2)$ but $F_\mu(\lambda, 0) = \lambda = F_\nu(\lambda, 0)$ for $\lambda \in \D$, it follows that $F_\mu \neq F_\nu$. Furthermore, the function $G := tF_\mu + (1-t)F_\nu$ is also a solution to the extremal problem \eqref{deg_red}. We shall show that $F \neq G$; then, nonuniqueness follows.

\par Suppose that $F = G$. According to Theorem \ref{kzw:1}, there exist $\tau_i, \omega_i$ for $i = 0, 1, 2$ (the modified problem does not require the application of automorphisms of the unit ball and unit disc) such that $F = F_{\tau_0, \omega_0}$ and $G = tF_{\tau_1, \omega_1} + (1-t)F_{\tau_2, \omega_2}$.
If $F = G$, then the set of singularities for both functions must be the same.

For any function $F_{\tau_i, \omega_i}$, where $i = 0, 1, 2$, the singularities lie on the parabola $z_2^2 = 2(1-\tau_i z_1)/\omega_i^2$. Since $F$, $F_\mu$, and $F_\nu$ each have uncountably many singularities, $F=tF_\mu+(1-t)F_\nu$, $F$ must share at least two singularities with at least one of $F_\mu$ or $F_\nu$. Without loss of generality, we assume that it shares two singularities with $F_\mu$, at the points $(z_1, z_2)$ and $(w_1, w_2)$, where $z_1 \neq w_1$, i.e., $F_\mu(z_1,z_2)=F(z_1,z_2)=\infty, F_\mu(w_1,w_2)=F(w_1,w_2)=\infty.$

From the parabola equations for the points above, we derive the following system of equations:
\[
\begin{cases}
    (1-\tau_1 z_1)/\omega_1^2=(1-\tau_0 z_1)/\omega_0^2\\
    (1-\tau_1 w_1)/\omega_1^2=(1-\tau_0 w_1)/\omega_0^2.
\end{cases}
\]
Dividing these equations side by side and multiplying the factors out, we obtain
$$(1-\tau_1z_1)(1-\tau_0w_1)=(1-\tau_1w_1)(1-\tau_0z_1).$$ 
Simplifying this expression, we get $(\tau_1-\tau_0)(z_1-w_1)=0,$ which means that $\tau_1=\tau_0.$
This, in turn, implies that $\omega_0^2 = \omega_1^2$. Therefore, $F = F_\mu$, and since $F = G$, it follows that $F = F_\mu = F_\nu$, which contradicts the assumption that $F_\mu \neq F_\nu.$

Hence, $F \neq G$, and the solution to the problem \eqref{deg_red} is not unique if the problem is degenerate, and such that $\omega\neq\omega_2$.

\subsection{Geodesics passing through points $(0, w_2), (1,0)$}
In the case of $\omega\neq\omega_2$, we will use complex 2-geodesics (in the further part of the article we will refer to them as complex geodesics for simplicity) that pass through points $(0,w_2)$ and $(1,0)$, where $w_2 \in (-1,1)$.

Recall that the complex geodesics in the unit Euclidean ball $\Bn{2}$ were described by Jarnicki, Pflug, and Zeinstra in \cite{9}:

\begin{thm}\label{jpz} Complex geodesics in 2-dimensional complex ellipsoid $\Bn{2}=\El(1,1)$ are affine discs. \end{thm}

A biholomorphism composed with a complex geodesic results in another complex geodesic within the target domain. Since the automorphisms of the unit ball $\Bn{2}$ are biholomorphisms, when they are composed with complex geodesics in $\Bn{2}$, the result remains a complex geodesic in $\Bn{2}$. One can show that any complex geodesic in $\Bn{2}$ is a composition of a horizontally flat geodesic $$\D\ni\lambda\mapsto(\lambda,0)\in\Bn{2}$$ with an automorphism of the unit ball $\Bn{2}.$

\par Using these facts, we construct the complex geodesic through points $(1,0)$ and $(0,w_2)$, where $w_2 > 0$. We begin with the geodesic $$\varphi_0:\D\ni\lambda\mapsto(\lambda,0)\in\Bn{2}.$$

Recall that the transformation \begin{align}\label{auto} \phi_a(z,w)=\left(\sqrt{1-a^2}\frac{z}{1+aw},\frac{w+a}{1+aw}\right),\ a\in(-1,1), \end{align} is an automorphism of the unit ball $\Bn{2}$. Composing $\phi_a$ and $\varphi_0$, one gets
$$(\phi_a\circ\varphi_0)(\lambda)=\left(\sqrt{1-a^2}\lambda, a\right),$$ which is a flat horizontal geodesic in $\Bn{2}.$

Consider the unitary automorphism generated by the matrix \[U_a:=
\begin{pmatrix}
\sqrt{1-a^2} & a \\
-a & \sqrt{1-a^2}
\end{pmatrix}
\] i.e., \begin{align}\label{uni} u_a(z,w):=\left(\sqrt{1-a^2}z+aw, -az+\sqrt{1-a^2}w\right). \end{align}

The composition \begin{align}\label{geo} \varphi_a(\lambda):=(u_a\circ\phi_a\circ\varphi_0 )(\lambda)=\left((1-a^2)\lambda+a^2, a\sqrt{1-a^2}(1-\lambda)\right) \end{align} is the desired complex geodesic.

Note that regardless of the choice of $a\in(0,1)$, the above geodesic always passes through the point $(1,0)$ - it suffices to set $\lambda=1.$ It remains to check that if $a$ is selected appropriately, the geodesic also passes through the point $(0,w_2).$

The first coordinate is zero if $\lambda = -\frac{a^2}{1-a^2}.$ For this choice of $\lambda$, the second coordinate is equal to $\frac{a}{\sqrt{1-a^2}}.$ Rearranging $\frac{a}{\sqrt{1-a^2}}=w_2$, we find that $a:=\frac{w_2}{\sqrt{1+w_2^2}}\in\D$ satisfies the requirement.

\subsection{Problem transformation}
Recall that we are considering the problem \begin{align}\label{problem2} \begin{cases} \Bn{2}\ni (0, 0) \mapsto 0\in\D, \\ \Bn{2}\ni (0, w_2) \mapsto \kappa\in\D, \\ \Bn{2}\ni (z_1, 0) \mapsto z_1\in\D. \end{cases}
\end{align} Using the rotations in $\Bn{2}$ and $\D$, we can assume that $w_2>0$.

Through the points $(0,0), (z_1,0)$ and $(1,0)$ there pass the geodesic $\varphi_0.$  Through $(0, w_2)$ and $(1,0)$ there pass a geodesic $\varphi_c,\ c:=\frac{w_2}{\sqrt{1+w_2^2}}.$ 

\par Our goal is to transform geodesics $\varphi_0$ and $\varphi_c$ to geodesics symmetric with respect to the $z$ axis, i.e, to transform $\varphi_0\mapsto\varphi_b,\ \varphi_c\mapsto\varphi_{-b},$ for some $b\in(0,1).$ This is possible because the family of automorphisms $(-1,1)\ni a\mapsto (u_a\circ\phi_a)(z,w),$ where $u_a$ is described by \eqref{uni} and $\phi_a$ by \eqref{auto}, is continuous, and the family fixes the point $(1,0).$


\par Below can be seen the above procedure: on the left, the initial problem and the geodesics passing through the interpolated points. On the right, the transformed problem - geodesics which "carried" interpolated points were transformed in such a way that those lie on "symmetrically" located geodesics $\varphi_{\pm b}.$

\begin{figure}[h!]
    \centering
    \includegraphics[width=0.8\textwidth]{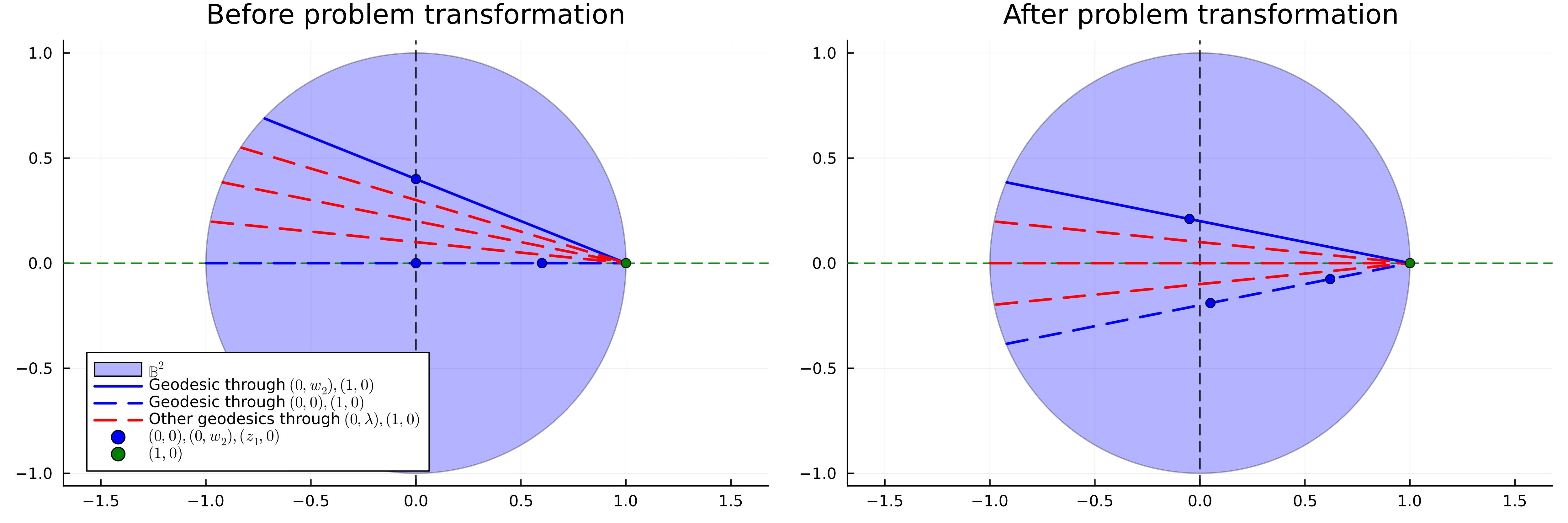}
\end{figure}

\subsection{Case: $\omega=\omega_2$}
In the previous subsection, using automorphisms of the unit ball, we transformed the Nevanlinna-Pick interpolation problem in such a way that the interpolated points lie on complex geodesics
$$\varphi_{\pm a}(\lambda)=\left((1-a^2)\lambda+a^2, \pm a\sqrt{1-a^2}(1-\lambda)\right)$$ for some $a\in(0,1).$ These geodesics cross at the point $(1,0).$ Additionally, they share the same left inverse from the class $\mathcal{F}_{D},$ namely $$F(z,w):=\frac{2z(1-z)-w^2}{2(1-z)-w^2}.$$  This left inverse is, in fact, also a left inverse for the geodesic $\varphi_0(\lambda)=(\lambda,0).$

  We will prove that there exists another extremal solution to the transformed problem such that it is not additionally a left inverse to $\varphi_0$.

Consider the complex ellipsoid $\El(1,1/2)$, where $$\mathcal{E}(p):=\left\{(z_1,...,z_d)\in\Cn{n}:\ \sum_{j=1}^{d}|z_j|^{2p_j}<1\right\}.$$  The corresponding geodesic to $\varphi_{\pm a}$ in that ellipsoid is $$\psi_{a}(\lambda)=\left((1-a^2)\lambda+a^2, a^2(1-a^2)(1-\lambda)^2\right).$$
 This follows as a consequence of the following "transport lemma" from \cite[Lemma 8]{9}:
\begin{lm}
    Let $\varphi = (\varphi_1, \dots, \varphi_d): \D \to \mathcal{E}(p)$ be a complex geodesic with $\varphi_j \not\equiv 0$, $j = 1, \dots, d$. Put
\[
B_j(\lambda) := e^{i \theta_j} \left( \frac{\lambda - \alpha_j}{1 - \overline{\alpha_j} \lambda} \right)^{k_j},\ \psi_j(\lambda)  :=e^{-i \theta_j} \left( a_j \frac{1 - \overline{\alpha_j} \lambda}{1 - \overline{\alpha_0} \lambda} \right)^{1/p_j}, \ \lambda\in\D,\ j = 1, \dots, d.
\]
Fix $\tilde{p} = (\tilde{p}_1, \dots, \tilde{p}_d) \in [1/2, +\infty)^d$ and define
\[
\tilde{\varphi}_j := B_j \psi_j^{p_j / \tilde{p}_j}, \quad (j = 1, \dots, d), \quad \tilde{\varphi} := (\tilde{\varphi}_1, \dots, \tilde{\varphi}_d).
\]
Then $\tilde{\varphi}$ is a complex geodesic in $\mathcal{E}(\tilde{p})$.

\end{lm}
To apply the above lemma and get the desired $\psi_a$, it suffices to take the parameters: $d=2, p_1=p_2=\widetilde{p_1}=1, \widetilde{p_2}=1/2,\theta_1=\theta_2=0, k_1=1,k_2=0, a_1=1-a^2, a_2=a\sqrt{1-a^2},\alpha_0=0,\alpha_1=\frac{a^2}{a^2-1},\alpha_2=1.$

Now we will construct a left inverse to $\psi_a$ using the Lempert theorem. Recall its statement \cite{99},  \cite[Section 8.2, Lemmas 8.2.2 and 8.2.4, Remark 8.2.3]{1}:

\begin{thm}(1982, Lempert)\label{lempert} If $\Omega$ is convex and bounded domain in $\Cn{d},$ and $\varphi\in\hols{\D}{\Omega},$ then $\varphi $ 
    is Kobayashi geodesic if and only if there exists a holomorphic function $F\in\hols{\Omega}{\D}$ such that $F\circ\varphi=id_{\D}.$ \par Additionally, if $\varphi$ is Kobayashi extremal, then there exists a function $h\in H^{1}(\D,  \mathbb{C}^d)$ such that 
    \begin{enumerate}
        \item[(i)] $\text{Re}[(z-\varphi(\lambda))\bullet(\overline{\lambda}h(\lambda))]<0,\ \text{for any}\ z\in\Omega, \ \text{and}\ a.e.\ \lambda\in\T, $
        \item[(ii)] a function $F$ can be found as a solution with respect to $\lambda$ of the equation $(z-\varphi(\lambda))\bullet h(\lambda)=0\ \text{for any } z\in\Omega,\ a.e.\ \lambda\in\T,$ 
    \end{enumerate}
     where $z\bullet w = \sum_{j=1}^{d}z_{j}w_{j},\ z,w\in\Cn{d}.$
\end{thm}
The ellipsoid $\El(1/2,1)$ is convex and bounded, so it is possible to define the function $h$ from the Lempert theorem as $$h(\lambda):=\lambda\rho(\lambda)\overline{\nu(\psi_a(\lambda))},\ \lambda\in\T,$$ where $\nu$ is the outer normal, and $\rho$ is a positive analytic function such that $h$ extends to an analytic function. This can be found in \cite[Remark 5]{9}.

In this case, $\nu(z,w):=(z,\frac{w}{2|w|}),\ (z,w)\in\partial\Bn{2}, w\neq0.$ One can compute that $\lambda\overline{\nu(\psi_a(\lambda))}= (a^2\lambda+1-a^2,-1/2).$ Therefore, it suffices to set $\rho\equiv1.$

The desired left inverse is a solution $\lambda=\lambda(z,w), \ (z,w)\in\El(1,1/2)$ to:
$$(a^2+(1-a^2)\lambda-z)(a^2\lambda+1-a^2)-1/2(a^2(1-a^2)(\lambda-1)^2-w))=0.$$ 
This solution is \begin{align*}
  \lambda(z, w) &= \frac{(1-a^2)^2-a^2(z-1)-\sqrt{[a^2(z-1)]^2+a^2(a^2-1)w+(1-a^2)^2}}{a^2(a^2-1)}.
\end{align*}
It is important to note that the square root is single-branched in $\El(1,1/2)$. For this we should see that any $a\in(0,1)$ if $$[a^2(z-1)]^2+a^2(a^2-1)w+(1-a^2)^2=0,$$ then $(z,w)\in\Cn{2}\setminus\El(1,1/2)$.
Computing out $w$ from the equation above we get parametrized form of the above equation 
$$\D\ni z\mapsto\left(z,\frac{a^4(1-z)^2+(1-a^2)^2}{a^2(1-a^2)}\right)\in\Cn{2}.$$
We have to show that for any $z\in\D$ and $a\in(0,1)$ one always has 
$$|z|^2+\frac{1}{a^2(1-a^2)}\left|a^4(1-z)^2+(1-a^2)^2\right|\geq1.$$
Setting $b=a^2,\ z=r(\cos\theta+i\sin\theta),\ r\in[0,1),\theta\in[0,2\pi),$ we get 
\begin{align*}r^2+\frac{1}{b(1-b)}(&((1-b)^2+b^2(1-2r\cos\theta+r^2(2\cos^2\theta-1))^2+\\&4b^4r^2\sin^2\theta(r\cos\theta-1)^2)^{1/2}.
\end{align*}
Substituting $u=\cos\theta, u\in[-1,1]$ problem reduces to finding a minimum of 
\begin{align*}  
r^2+\frac{1}{b(1-b)}(&((1-b)^2+b^2(1-2ru+r^2(2u^2-1))^2\\ &+4b^4r^2(1-u^2))(ru-1)^2)^{1/2}
\end{align*} with constraints
$r\in[0,1),u\in[-1,1], b\in(0,1).$ It can be checked that this function does not have extreme points inside the cube $(0,1)\times(-1,1)\times(0,1).$ It attains its minimum on the boundary for $(r,u)=(1,1), b=1$ and it is equal to $1.$ Which is enough to conclude that the square root is analytic in $\El(1,1/2).$

Finally, it suffices to observe that $\widetilde{F}(z,w):=\lambda(z,w^2),\ \widetilde{F}\in\hols{\Bn{2}}{\D}$ is a left inverse to both $\varphi_{\pm a}.$ However, $\widetilde{F}$ is not a left inverse to the geodesic $\varphi_0:\zeta\mapsto(\zeta,0),$ while $F$ is its left inverse. Therefore, $\widetilde{F}$ generates a different solution from $F$, and nonuniqueness follows.

\section{Acknowledgments}  The author is grateful to Włodzimierz Zwonek from Jagiellonian University in Cracow for sharing his ideas and help. The author would like thank also to Pascal Thomas from Institut de Mathématiques de Toulouse, France, for his review and his precious and detailed remarks.

\address{Faculty of Mathematics and Computer Science, Jagiellonian University, Lojasiewicza 6, 30-348, Krakow, Poland}
\email{dariusz.piekarz@doctoral.uj.edu.pl}

\end{document}